\documentclass[12pt,english,reqno]{amsart}
\usepackage{amssymb,amsmath,amstext,amsthm,amsfonts}
\usepackage[dvips]{graphicx}
\usepackage[ansinew]{inputenc}

\newtheorem{Theorem}{Theorem}

\newtheorem*{Maintheorem}{Main Theorem}

\newtheorem{maintheorem}{Theorem}

\newcommand{\cmt}{\begin{maintheorem}}
\newcommand{\fmt}{\end{maintheorem}}

\newtheorem{maincorollary}[maintheorem]{Corollary}

\newcommand{\cmc}{\begin{maincorollary}}
\newcommand{\fmc}{\end{maincorollary}}

\newtheorem{T}{Theorem}[section]
\newcommand{\cte}{\begin{T}}
\newcommand{\fte}{\end{T}}

\newtheorem{Corollary}[T]{Corollary}
\newcommand{\cco}{\begin{Corollary}}
\newcommand{\fco}{\end{Corollary}}

\newtheorem{Proposition}[T]{Proposition}
\newcommand{\cpr}{\begin{Proposition}}
\newcommand{\fpr}{\end{Proposition}}

\newtheorem{Lemma}[T]{Lemma}
\newcommand{\cle}{\begin{Lemma}}
\newcommand{\fle}{\end{Lemma}}

\newcommand{\csle}{\begin{Lemma}}
\newcommand{\fsle}{\end{Lemma}}

\newtheorem{Remark}[T]{Remark}
\newcommand{\cre}{\begin{Remark}}
\newcommand{\fre}{\end{Remark}}

\newtheorem{Definition}[T]{Definition}
\newcommand{\cde}{\begin{Definition}}
\newcommand{\fde}{\end{Definition}}

\def \BB {{\mathbb B}}
\def \AA {{\mathbb A}}
\def \RR {{\mathbb R}}

\def \NN {{\mathbb N}}

 \def \cp {\mathcal{P}}
 \def \cc {\mathcal{C}}
 
 \def \cs {\mathcal{S}}
 \def \cu {\mathcal{U}}

\newcommand{\dem}{\begin{proof}}
\newcommand{\cqd}{\end{proof}}

\newcommand{\leb}{\operatorname{Leb}}
\newcommand{\adj}{\operatorname{adj}}
\newcommand{\dist}{\operatorname{dist}}
\newcommand{\supp}{\operatorname{supp}}

\begin{document}

\author{Vilton  Pinheiro}
\address{Departamento de Matem\'atica, Universidade Federal da Bahia\\
Av. Ademar de Barros s/n, 40170-110 Salvador, Brazil.}
\email{viltonj@ufba.br}

\thanks{Work carried out at the  Federal University of
Bahia (Brazil) and ICTP (Italy)
Partially supported by PADCT/CNPq(Brazil)}

\date{\today}

\title{SRB measures for weakly expanding maps}

\maketitle

\begin{abstract}
We construct SRB measures for endomorphisms satisfying conditions
far weaker than the usual non-uniform expansion. As a consequence,
the definition of a non-uniformly expanding map can be weakened. We
also prove the existence of an absolutely continuous invariant
measure for local diffeomorphisms, only assuming the existence of
hyperbolic times for Lebesgue almost every point in the manifold.
\end{abstract}

\section{Introduction}
In \cite{Ke}, Keller proved the existence of an absolutely
continuous invariant measure for any {\em non-flat multimodal map}
$f$, whenever it satisfies the {\em negative Schwarzian derivative
condition} and there exists a constant $\lambda$ such that for
Lebesgue almost all $x$
\begin{equation}\label{e1}
\limsup_{n\to\infty}\frac{1}{n}\sum_{j=0}^{n-1}\log|Df(f^j(x))|\ge\lambda>0.
\end{equation}

Besides that, Keller showed the existence of a finite number of
these measures whose union of  basins contains Lebesgue almost all
point of the domain of $f$. This measure is called a {\em physical}
or {\em SRB} (Sinai-Ruelle-Bowen) measure.

The result of Keller was somewhat generalized by
Alves-Bonatti-Viana~\cite{ABV} in the context of non-flat maps
defined on a compact Riemannian manifold $M$. For this, the
distortion control condition of {\em negative Schwarzian derivative}
was replaced by the so called {\em slow recurrence to the critical
set} and it was assumed that
\begin{equation}\label{e2}
\liminf_{n\to\infty}\frac{1}{n}\sum_{j=0}^{n-1}\log(\|(Df(f^j(x)))^{-1}\|^{-1})\ge\lambda>0
\end{equation}
for Lebesgue almost all $x\in M$. Of course, if the dimension of $M$
is one we have $|Df(x)|=\|(Df(x))^{-1}\|^{-1}$. Observe that
$\|(Df(x))^{-1}\|^{-1}>1$ means that $Df(x)$ expands in all
directions, that is, $\|(Df(x))v\|\ge\gamma\|v\|$ $\forall\,v\in T_x
M$ (where $\gamma=\|(Df(x))^{-1}\|^{-1}$). Systems
satisfying~(\ref{e2}) are called {\em non-uniformly expanding}, they
generalize the {\em uniformly expanding systems}. As particular
examples of this kind of systems we can mention one-dimensional maps
with {\em positive Lyapunov exponents} (like {\em quadratic maps}
and, in general, {\em non-flat multimodal maps}~\cite{MvS}), and in
higher dimension, the {\em Viana maps}~\cite{V}.

We note that any multimodal map that satisfies Keller's hypothesis
also satisfies {\em the slow recurrence condition to the critical
set}. That is, if $f$ is a $C^3$ non-flat multimodal map with {\em
positive Lyapunov exponent} (i.e., it satisfies (\ref{e1})) then the
hypothesis of {\em negative Schwarzian derivative} is stronger than
the condition of {\em slow recurrence to the critical set} (see
Proposition~\ref{Schwarzian}).

The purpose of this paper is to generalize the result of Keller  in
the context of \cite{ABV} (replacing the {\em negative Schwarzian
derivative} by the {\em slow recurrence to the critical set}) with
the far more weak condition of $\limsup$ as it appears in Keller's
theorem, i.e.,
\begin{equation}\label{e3}
\limsup_{n\to\infty}\frac{1}{n}\sum_{j=0}^{n-1}\log(\|(Df(f^j(x)))^{-1}\|^{-1})\ge\lambda>0.
\end{equation}
Moreover, we prove that the condition~(\ref{e3}) indeed implies
(\ref{e2}) and so, the definition of {\em non-uniformly expanding
map} can be weakened.

\subsection{Statement of results} Let $M$ be a compact Riemannian manifold
of dimension \( d\geq 1\) and \( \leb \) a normalized Riemannian volume
form on $M$ that we call {\em Lebesgue measure}.

A map $f\colon M\to M$  will be called non-flat if $f$ is local
$C^2$ diffeomorphism in the whole manifold  except  in a {\em
non-degenerate critical (or singular) set}  $\mathcal{C}\subset M$.
The definition of {\em non-degenerate critical set} is given at the
beginning of section~\ref{Hyperbolic Times} (If $dim(M)=1$ and $f$
satisfies the usual definition of non-flatness (see \cite{MvS}),
then it also satisfies the definition given above).

We say that $f$ satisfies the {\em
slow approximation condition}  if given any $\epsilon>0$ there exists $\delta>0$
such that for Lebesgue almost every point $x\in M$ we have
\begin{equation} \label{e.faraway1}
\limsup_{n\to+\infty}
\frac{1}{n} \sum_{j=0}^{n-1}-\log \mbox{dist}_\delta(f^j(x),\cc)
\le\epsilon,
\end{equation}
where \(\dist_{\delta}(x,\cc) \) denotes the \( \delta
\)-\emph{truncated} distance from \( x \) to \( \cc \) defined as \(
\dist_{\delta}(x,\cc) =  \dist(x,\cc) \), if \( \dist(x,\cc) \leq
\delta\), and \( \dist_{\delta}(x,\cc) =1 \), otherwise.

We call the {\em basin} of some invariant measure $\nu$ the set
$\mathcal{B}(\nu)$ of the points $x\in M$ such that the average of
Dirac measures along the orbit of $x$ converges in the weak$^*$
topology to $\nu$, that is,
$$\lim_{n\to +\infty}\frac{1}{n}\sum_{j=0}^{n-1}\phi(f^j(x))=
\int\phi\,d\nu\mbox{, }\forall\phi\in C^0(M).$$

\begin{Theorem}\label{SRB}
Let $f:M\to M$ be a non-flat  map satisfying the slow
approximation condition. If $f$ ( or some fixed iterate ) satisfies
\begin{equation}\label{e4}
\liminf_{n\to\infty}\frac{1}{n}\sum_{i=0}^{n-1}\log(
\|(Df(f^{i}(x)))^{-1}\|)\le-\lambda<0
\end{equation}
(or equivalently satisfies~(\ref{e3})) for Lebesgue almost every
point $x\in M$, then there exists a finite collection of ergodic
absolutely continuous invariant measures such that almost every
point in $M$ belongs to the basin of one of these measures.
\end{Theorem}

Theorem~\ref{SRB} is in fact a consequence of the Main Theorem below
that guarantees the existence of a global Markov structure with
integrable time function. A global Markov structure is composed by a
time function $x\mapsto R(x)\in\NN\cup\{\infty\}$, an induced map
$F(x)=f^{R(x)}(x)$ defined almost everywhere and a countable
partition  refining a finite triangulation of $M$. Each element of
this partition is sent by $F$, with good properties (see
section~\ref{Markov Structures} for details),  onto an element of
the triangulation. One can prove that such $F$ has an absolutely
continuous invariant measure $\nu$. Moreover, whenever $R$ is
$\nu$-integrable, it generates an absolutely continuous
$f$-invariant finite measure. We observe that the existence of a
Markov structure allows a more deep study of the dynamical
properties of the map $f$. For instance, it was used in
\cite{ALP,BY2,G,Y1,Y2} to study the decay of correlations and prove
the Central Limit Theorem for a large class of maps.

\begin{Maintheorem}
Every map satisfying the hypothesis of theorem~\ref{SRB} has a global
Markov Structure with integrable time function.
\end{Maintheorem}

The proof of the main theorem is a mix of the strategy adopted by
Alves-Luzzatto-Pinheiro in \cite{ALP} for  non-uniform expanding
maps (inspired in Young's paper \cite{Y1}) and the proof of Keller's
theorem in \cite{MvS}. Here we are able to simplify the proof of the
existence of the Markov Structure that appears in \cite{ALP,G,Y1}
and make it closer to the one dimensional case.

A map $f\colon M\to M$ is called {\em non-uniformly expanding} if
$f$ is a non-flat map satisfying the slow approximation condition
and it (or some fixed iterate) satisfies, for Lebesgue almost every
point $x\in M$, the following condition,
\begin{equation}\label{e5}
\limsup_{n\to\infty}\frac{1}{n}\sum_{i=0}^{n-1}\log(
\|(Df(f^{i}(x)))^{-1}\|)\le-\lambda<0
\end{equation}
(or equivalently satisfies~(\ref{e2})).

\begin{Theorem}\label{weaklyXnue}
A $C^2$ endomorphism  is non-uniformly expanding if and only if it
satisfies the hypothesis of theorem~\ref{SRB}.
\end{Theorem}

We want to remark that Theorem~\ref{weaklyXnue} deals with a more
restricted class of maps than Theorem~\ref{SRB}. For instance, maps
with a singular set $\cs$ ($|\det Df|_{\cs}|=+\infty$) are allowed
in the hypothesis of Theorem~\ref{SRB}. This difference happens
because, in the proof of Theorem~\ref{weaklyXnue}, we need a $C^2$
map to guarantee the integrability of $\log\|(Df)^{-1}\|$ with
respect to any absolutely continuous invariant measure (see
Lemma~\ref{PDL}).

A crucial ingredient in the proof of these results is the existence,
for almost all point $x\in M$, of moments $n=n(x)\in\NN$ such that
$f$ looks like an uniformly expanding map on some neighborhood
$V_n(x)$ of $x$, that is, this  neighborhood is sent by $f^n$, with
good properties of expansion and distortion, to some big ball
$B(f^n(x))$ centered at $f^n(x)$. The last result of this paper is
about hyperbolic times which are associated with these moments of
expansion mentioned above (see Proposition~\ref{Prop3}). In
\cite{AA} Alves-Ara\'ujo proved the existence of an {\em SRB}
measure when the critical set is empty and the first hyperbolic time
is Lebesgue integrable. Here, we were able to remove  the hypothesis
of integrability.

\begin{Theorem}\label{HTSRB}
Let $f:M\to M$ be a $C^2$ covering map (local diffeomorphism) on a
compact manifold M. If the first hyperbolic time function is defined
for  Lebesgue almost every point of M, then $f$ is a non-uniformly
expanding map.
\end{Theorem}

\section{Hyperbolic Times}
\label{Hyperbolic Times}

\begin{Definition}\label{criset}
Let $f:M\to M$ be a $C^2$  local diffeomorphism in the whole
manifold except  in  a  critical (or singular) set
$\mathcal{C}\subset M$. We say that $\mathcal{C}\subset M$ is a {\em
non-degenerate critical set}, more precisely, a $\beta$-{\em
non-degenerate critical set} ($\beta>0$) if $\exists B>0$ such that
the following three conditions hold:

\begin{enumerate}
 \item[(1)]
\quad $\displaystyle{\frac{1}{B}\dist(x,\mathcal{C})^{\beta}\leq
\frac {\|Df(x)v\|}{\|v\|}\leq B\, \dist(x,\mathcal{C})^{-\beta}}$
for all $x\in M\setminus\mathcal{C}$ and $v\in T_x M$.
\end{enumerate}

For every $x,y\in M\setminus\mathcal{C}$ with
$\dist(x,y)<\dist(x,\mathcal{C})/2$ we have
\begin{enumerate}
\item[(2)] \quad $\displaystyle{\left|\log\|Df(x)^{-1}\|-
\log\|Df(y)^{-1}\|\:\right|\leq
\frac{B}{\dist(x,\mathcal{C})^{\beta}}\dist(x,y)}$;
 \item[(3)]
\quad $\displaystyle{\left|\log|\det Df(x)|- \log|\det
Df(y)|\:\right|\leq
\frac{B}{\dist(x,\mathcal{C})^{\beta}}\dist(x,y)}$.
 \end{enumerate}
\end{Definition}

\begin{Proposition}\label{Schwarzian}
Let $f:I\to I$ be a non-flat $C^3$ multimodal map of the interval
$I$ with positive Lyapunov exponent, i.e., it satisfies (\ref{e1})
for almost every point $x\in I$. If $f$ satisfies the negative
Schwarzian derivative condition, then it also satisfies the
condition of slow recurrence to the critical set.
\end{Proposition}

\dem By Keller's result, there exists a finite collection
$\mu_1,\ldots,\mu_s$ of ergodic absolutely continuous invariant
measures such that $Leb(I\setminus\mathcal
B(\mu_1)\cup\ldots\cup\mathcal B(\mu_s))=0$. Moreover, the support
of each $\mu_j$ is a finite union of intervals and for Lebesgue
almost all $p\in\mathcal B(\mu_j)$ there is some $n_0=n_0(p)$ such
that $q=f^{n_0}(p)$ is a $\mu_j\,$-generic point, in particular
$q\in\supp\mu_j$. As $f$ is $C^3$, there are constants $A_0$,
$A_1>0$ such that $|Df(x)|<A_0\dist(x,\cc)<A_1$, $\forall\,x\in I$.
On the other hand, as $\log|Df|\in L^1(\mu_j)$ (Remark 1.2 of
\cite{Liu}), it follows that $\varphi:x\mapsto\log\dist(x,\cc)$ is
$\mu_j\,$-integrable.


Taking $I_{\delta}=\{x\in I\|\dist(x,\cc)<\delta\}$, we have
$\leb(I_{\delta})\searrow 0$. So, $\mu_j(I_{\delta})\to 0$ and also
$\lim_{\delta\to 0}\int\varphi_{\delta}d\mu_j=\lim_{\delta\to
0}\int_{I_{\delta}}\varphi d\mu_j=0$, where
$\varphi_{\delta}=\log\dist_{\delta}(\cdot\,,\cc)$. Thus, given
$\varepsilon>0$ let $\delta>0$ be such that
$-\int\varphi_{\delta}d\mu_j$ $\le\varepsilon$,
$\forall\,j=1,\ldots,s$. By Birkhoff's ergodic theorem,
$$\lim\frac{1}{n}\sum_{j=0}^{n-1}-\log\dist_{\delta}(f^j(q),\cc)=
-\int\varphi_{\delta}d\mu_j\le\varepsilon,$$ for any
$\mu_j\,$-generic $q$. Thus,
$$\lim\frac{1}{n}\sum_{j=0}^{n-1}-\log\dist_{\delta}(f^j(p),\cc)=$$
$$=\lim\frac{1}{n}\sum_{j=0}^{n-1}
-\log\dist_{\delta}(f^j(f^{n_0(p)}(p)),\cc)\le\varepsilon,$$ for
Lebesgue almost all $p\in\mathcal B(\mu_1)\cup\ldots\cup\mathcal
B(\mu_s)$. \cqd

Let us fix $0<b<\frac{1}{2}\min\{1,1/\beta\}$. If $f:M\to M$ is a
$C^2$ local diffeomorphism outside a $\beta$-non-degenerate critical
set $\mathcal{C}$ then, given $0<\sigma<1$ and $\delta>0$, we will
say that $n$ is a $(\sigma,\delta)$-{\em hyperbolic time} for a
point $x\in M$ if for all $1\le k\le n$ we have
$\prod_{j=n-k}^{n-1}\|(Df\circ f^j(x))^{-1}\|\le {\sigma}^k\mbox{
and }$ $\dist_{\delta}(f^{n-k}(x),\mathcal{C})\ge {\sigma}^{b k}$.
We denote the set of points of $M$ such that $n\in\NN$ is
$(\sigma,\delta)$-{\em hyperbolic time} by $H_n(\sigma,\delta)$.

\begin{Proposition}\cite{ABV}
\label{Prop2} Let $f:M\to M$ be a $C^2$ non-flat map satisfying the
slow approximation condition. Given $\lambda>0$ there exist
$\delta>0$ and $\theta>0$, depending only on $f$ and $\lambda$, such
that $$\#\{1\le j\le n\ \|\ x\in
H_j(e^{-\lambda/4},\delta)\}\ge\theta\,n,$$ whenever
$\sum_{i=0}^{n-1}
    \log \|(Df(f^{i}(x)))^{-1}\|^{-1}>\lambda\,n$.
\end{Proposition}

\begin{Corollary}
\label{Cor1} Let $f:M\to M$ be a $C^2$ non-flat map satisfying the
slow approximation condition. Given $\lambda>0$ there exist
$\delta>0$ and $\theta>0$, depending only on $f$ and $\lambda$, such
that if $A\subset M$ with $\limsup\frac{1}{n}\sum_{i=0}^{n-1} \log
\|(Df(f^{i}(x)))^{-1}\|^{-1}>\lambda$ for Lebesgue almost all $x\in
A$ then
$$\limsup_n\frac{1}{n}\#\{1\le j\le n\ \|\ x\in
 H_j(e^{-\lambda/4},\delta)\}\ge\theta$$ for Lebesgue almost all
 $x\in A$
\end{Corollary}

We finish this section stating a proposition that assures a good
behavior, with respect to $f^n$, of a neighborhood of a point $x$
when $n$ is a hyperbolic time for this point.

\begin{Proposition}\label{Prop3}
Let $f:M\to M$ be a $C^2$ local diffeomorphism outside a
non-degenerate critical set $\mathcal{C}$. Given $\sigma<1$ and
$\delta>0$, there exist $\delta_1,\rho
>0$, depending only on $\sigma,\delta$ and on the map $f$, such
that for any $x\in M$ and \( n\geq 1 \) a
$(\sigma,\delta)$-hyperbolic time for \( x \), there exists a
neighborhood \( V_n(x) \) of \( x \) with the following
properties:
 \begin{enumerate}
 \item $f^{n}$ maps $V_n(x)$ diffeomorphically onto the ball
  $B_{\delta_1}(f^{n}(x))$;
 \item $\dist(f^{n-j}(p),f^{n-j}(q)) \le
  \sigma^{j/2}\dist(f^{n}(p),f^{n}(q))$ $\forall p, q\in V_n(x)$ and $1\le j<n$;
 \item $\log |\frac{\det Df^n(p)}{\det Df^n(q)}|\le \rho\,\dist(f^n(p),f^n(q))$
  $\forall\,p$, $q\in V_n(x)$.
\end{enumerate}
\end{Proposition}

We shall often refer to the sets \( V_n(x) \) as \emph{hyperbolic
pre-balls} and to their images \( f^{n}(V_n(x)) \) as
\emph{hyperbolic balls}.  Notice that the latter are indeed balls
of radius \( \delta_1 \).

\dem For the proofs of items 1 and 2 see Lemma 5.2 in \cite{ABV}.
Let $p$, $q\in V_n(x)$. As $n$ is a hyperbolic time for $x$, it
follows from item (2) above that $\dist(f^j(q),\mathcal{C})\ge
\dist(f^j(x),\mathcal{C})- \dist(f^j(x),f^j(q))\ge
\sigma^{b(n-j)}-\delta_1\sigma^{(n-j)/2}\ge
(1-\delta_1)\sigma^{b(n-j)}$. Now, using condition (3) of the
definition of a non-degenerate critical set, we get
$$\log \bigg|\frac{\det Df^n(q)}{\det Df^n(p)}\bigg|\le
\sum_{j=0}^{n-1}\log \bigg|\frac{\det Df(f^{j}(q))}{\det Df(f^{j}(p))}\bigg|\le$$
$$\le \sum_{j=0}^{n-1}B\,\frac{\sigma^{(n-j)/2}\dist(f^n(q),f^n(p))}{((1-\delta_1)
\sigma^{b(n-j)})^\beta}\le$$
$$
\le\frac{B}{(1-\delta_1)^\beta(1-\sigma^{1/2-b\beta})} \dist(f^n(q),f^n(p))
$$

\cqd

\section{Markov Structures and The Partitioning Algorithm}

In this section we will say what we mean by a {\em Global Markov
Structure} and  prove the first part of the {\bf Main Theorem}. A
{\bf Global Markov Structure} for a map $f:M\to M$ is an induced
{\em piecewise uniformly expanding Markovian map\/} defined in an
open subset of $M$ with full Lebesgue measure. Precisely, we will
say that $f$ has a Global Markov Structure if there exists an open
subset $M'\subset M$ with $\leb(M\setminus M')=0$ and a function
$R:M'\to\NN\cup\{+\infty\}$ (called {\bf time function}) such that
the induced map $F:M'\to M$ given by $F(x)=f^{R(x)}(x)$ is a {\em
piecewise uniformly expanding Markovian map\/} defined on $M$ (see
 the next paragraph). The pair $(F,R)$ is
called a Global Markov Structure for $f$.

\subsection{Piecewise expanding Markovian map}
{\em A $C^2$ piecewise uniformly expanding Markovian map\/} defined
on $M$ is a map $F$ defined in an open subset of $M$ with full
Lebesgue measure such that there is a (mod $0$) countable partition
$\cp'$ refining a finite partition (mod $0$)
$\cp=\{P_1,\ldots,P_n\}$, such that each element is given as the
interior of some element in a triangulation $\mathcal T$ of $M$,
with $F$, $\cp'$ and $\cp$ satisfying the following properties:
\begin{enumerate}
\item $P_j$ has a piecewise $C^2$ boundary with finite
$(d-1)$-dimensional volume,  $\forall\,j=1,\ldots,n$.
\item $\exists\,0<\kappa<1$ such that $\|D F(x)^{-1}\| <\kappa$,
$\forall\,x\in U$ and $\forall\,U\in\cp'$.
\item $\forall\,U\in\cp'$, $F\vert_U$ is a $ C^{2} $ diffeomorphism  onto some element
of $\cp$.
\item $\exists\,K>0$ such that
$ \log\left|\frac{\det DF(x)}{\det DF(y)}\right| \leq K \dist(F(x),
F(y))$, $\forall\,x, y\in U$ and $\forall\,U\in\cp'$.
\end{enumerate}

The theorem below assures that every $C^2$ piecewise uniformly
expanding Markovian map $F\colon\Delta\to\Delta$ has an absolutely
continuous invariant measure $\nu$ whose density belongs to
$L^{\infty}(Leb)$; see e.g. Lemma $4.4.1$ of \cite{Aa}. Moreover, it
is straightforward to check that if $R$ is $\nu$-integrable, then
$$ \mu =\sum_{j=0}^{\infty}f_{\ast}^j\left(\nu\vert \{R>j\}\right)$$
 is an absolutely continuous $f$-invariant measure. Here  $\nu\vert
\{R>j\}$  denotes the measure given by $\nu\vert
\{R>j\}(A)=\nu(A\cap\{R>j\})$, and $f_*^j$ denotes the
push-forward of the measure by $f^j$.

We say that the {\bf Markov Structure has integrable time function}
if $R$ is integrable with respect to any absolutely continuous
$F$-invariant measure. As a consequence to the theorem below, if we
want to show that a Markov Structure has integrable time function,
we only need to verify the integrability of $R$ with respect to the
finite collection of ergodic absolutely continuous invariant
measures given by this theorem.

\begin{Theorem} \cite{Aa,Al,AV}\label{MST}
If $F\colon\Delta\to\Delta$ is a $C^2$ piecewise uniformly expanding
Markovian map, then there exists a finite set of ergodic absolutely
continuous invariant measures such that Lebesgue almost every point
in $\Delta$ belongs to the basin of one of these measures. Moreover,
the density of each of these measures with respect to Lebesgue is
uniformly bounded by some constant.
\end{Theorem}

\subsection{The Partitioning Algorithm}

Let $f:M\to M$ be a $C^2$ non-flat  map satisfying the slow
approximation condition and take $\lambda_0>0$. Suppose that
$$\limsup_{n\to\infty}\frac{1}{n}\sum_{i=0}^{n-1} \log
\|(Df(f^{i}(x)))^{-1}\|^{-1}>\lambda_0$$ for Lebesgue almost every
$x\in M$. In this section we will show how one can construct a
Markovian Structure on $M$.

Replacing, if is necessary, $f$ by some $f^{n_0}$
we may assume, without loss of generality, that $\lambda_0> 16\log2$.

Let $\lambda=\lambda_0/4$ and
\begin{equation}\label{sigma}
\sigma=\exp(-\lambda)<\frac{1}{16}
\end{equation} be fixed from now on.
Let $\delta>0$ be given by Proposition~\ref{Prop2} and let
$\delta_1>0$ be the radius of the hyperbolic ball given by
Proposition~\ref{Prop3}. For short, let $H_j=H_j(\sigma,\delta)$,
$\forall\,j\ge 1$.
\begin{Remark}
\label{remarkTHETA} It follows from Corollary~\ref{Cor1} the
existence of a constant $\theta>0$ such that
$\limsup_n\frac{1}{n}\#\{1\le j\le n\ \|\ x\in
 H_j\}\ge\theta$ for Lebesgue almost all
 $x\in M$.
\end{Remark}

Choose some $0<\delta_0<\delta_1/4$ and a finite partition $\cp$ of
$M$ (mod 0) generated by the interior of the elements of a
triangulation of $M$ with diameter smaller than $\delta_0$, that is,
$diameter(P)<\delta_{0}$ $\forall P\in\cp$.

Given $U\subset M$ and $0<r$ define the $r$-neighborhood of $U$ as
$B_r(U)=\{x\in M\,\|\,\dist(x,U)<r\}$. For each element $P\in\cp$,
set $P^{1}=B_{\delta_{0}}(P)$ and let $\{I^1_k(P)\}_{k\in\NN}$ be
the partition (mod 0) of $P^{1}\setminus P$ into the collection of
``rings''
$$I^1_{k}(P)=\left\{x\in P^{1}\setminus P\,\|\,\delta_{0}\,\sigma^{k/2} \le
\dist(x,\partial P) < \delta_{0}\,\sigma^{(k-1)/2}\right\},$$ where
$k=1,2,3,...$.

Now we present an algorithm to construct a refinement $\mathcal{P}'$
of $\mathcal{P}$ such that each element $Q\in\mathcal{P}'$ will be
sent, at a hyperbolic time, diffeomorphically onto some element of
$\mathcal{P}$. That is, each $Q\in\mathcal{P}'$ will be contained in
some hyperbolic pre-ball $V_n(q)$ for some $q\in Q\cap H_n$ and
there will be some $P\in\mathcal P$ such that
$Q=(f^n|V_n(q))^{-1}(P)$. The major difficulty of the algorithm is
to prevent overlaps of the new elements with the elements
constructed in previous steps. To deal with this problem, we will
define inductively, for each step $n\in\NN$, an auxiliary function
$t_n$ (we are assuming $0\in\NN$). This function will tell us if it
is possible to use, at the step $n$, a point $x\in H_n$ as the
``center'' of a hyperbolic pre-ball $V_n(x)$ to generate a new
element of the partition.

Given an element $P_0\in\cp$, let us set $\Delta_0(P_0)=P_0$. We
define a function $t_0:P_0\to\NN$ by
$$t_0(x)=\min\{k\in\NN\,\|\,\delta_0
\sigma^{k/2}\le\dist(x,\partial P_0)\},$$ and the subsets
$\AA_0(P_0)$ and $\BB_0(P_0)$ of $\Delta_0(P_0)$ given by
$$\AA_0(P_0)= \{x\in\Delta_0(P_0)\,\|\, t_0(x) = 0\}$$
$$\mbox{and}$$
$$\BB_0(P_0)=
\{x\in\Delta_0(P_0)\,\|\, t_0(x) > 0\}.$$ In fact, as
$diameter(P_0)<\delta_0$, we have $\AA_0(P_0)=\emptyset$ and
$\BB_0(P_0)=P_0$.

We define inductively  the sets $\Delta_{n}(P_0)$, $\AA_{n}(P_0)$,
$\BB_{n}(P_0)$ and functions $t_{n}:P_0\to\NN$ for every $n\in\NN$.
For this, let us assume  that for some $n>0$ the sets
$\Delta_{k}(P_0)$, $\AA_{k}(P_0)$, $\BB_{k}(P_0)$ and functions $
t_{k}:P_0\to\NN$ are already defined for all $0\le k\le n-1$.

Our construction will depend on how many points of $\AA_{n-1}(P_0)$
have $n$ as a hyperbolic time. So, for the definition of
$\Delta_n(P_0)$ and $t_n$, we consider two cases.
\subsection*{First case} If
$\leb(\AA_{n-1}(P_0)\cap H_n)=0$ we set
$\Delta_n(P_0)=\Delta_{n-1}(P_0)$. Setting $t'_n:P_0\to\NN$ as
$t'_n(x)=0$ $\forall\,x$, we define $t_{n}:P_0\to\NN$ as
$$t_n(x)=\max\{t'_n(x),t_{n-1}(x)-1\}\mbox{
}(=\max\{0,t_{n-1}(x)-1\}).$$ Thus, we may write
\begin{equation*} t_{n}(x)=
\begin{cases}
0 & \text{ if } x\in \AA_{n-1}(P_0)\\
t_{n-1}(x)-1 & \text{ if } x\in \BB_{n-1}(P_0)
\end{cases}.
\end{equation*}

\subsection*{Second case}
If $\leb(\AA_{n-1}(P_0)\cap H_n)>0$, then, for Lebesgue almost every
point  $x\in \AA_{n-1}(P_0)\cap H_n$, there exists some $P\in\cp$
such that $f^n(x)\in P$. As the diameter of $P^1$ is smaller than
$3\delta_1/4$ we also have $B_{\delta_1}(f^n(x))\supset P^1\supset
P$. Proposition~\ref{Prop3} assures that there exists a hyperbolic
pre-ball $V_n(x)$ such that $f^n|V_n(x):V_n(x)\to
B_{\delta_1}(f^n(x))$ is a diffeomorphism with uniform bounded
distortion ( not depending on $x$ or $n$). Denote by
$\mathcal{U}_n(P_0)$ the collection of all sets
$(f^n|V_n(x))^{-1}(P)$ contained in $P_0$, where
$x\in\AA_{n-1}(P_0)\cap H_n$ and $f^n(x)\in P\in\cp$.

{\em We claim that the elements of $\mathcal{U}_n(P_0)$ are mutually
disjoint.} To check this, let us assume that $x_1$, $x_2\in H_n\cap
P_0$, $f^n(x_1)\in P_1\in\cp$ and $f^n(x_2)\in P_2\in\cp$. Also, let
$U=(f^n|_{V_n(x_1)})^{-1}(P_1)$ and $V=(f^m|_{V_m(x_2)})^{-1}(P_2)$.
As the elements of $\cp$ are small compared to
$B_{\delta_1}(f^n(x_1))$ or $B_{\delta_1}(f^n(x_2))$
($diameter(P)<\delta_0<\delta_1/4$ $\forall\,P\in\cp$), if
$P_1=P_2$, then
$(f^n|_{V_n(x_1)})^{-1}|_{P_1}=(f^n|_{V_n(x_2)})^{-1}|_{P_2}$ and
so, $U=V$. On the other hand, if $P_0\ne P_1$, then $f^n(U)\cap
f^n(V)=P_0\cap P_1=\emptyset$ and so, $U\cap V=\emptyset$, proving
the claim.

Given $U\in\mathcal{U}_n(P_0)$, with $U=(f^n|V_n(x))^{-1}(P)$ and
$P\in\cp$, let $$U^1=(f^n|V_n(x))^{-1}(P^1)$$ and, for
$j=1,2,3,...$, $$U^1(j)=(f^n|V_n(x))^{-1}(I^1_j(P)).$$

Define
$$\Delta_{n}(P_0)=\Delta_{n-1}(P_0)\setminus
\bigcup_{U\in\mathcal{U}_n(P_0)} U.$$

For each point $x\in\bigcup_{U\in\mathcal{U}_n(P_0)}U^1\setminus U$,
set
$$t'_n(x)=\max\{i\,\|\, x\in U^1(i)\mbox{ and
}U\in\mathcal{U}_n(P_0)\},$$ and for the other points $x\in P_0$,
set $t'_n(x)=0$.

As before, define $t_{n}:P_0\to\NN$ by
$$t_n(x)=\max\{t'_n(x),t_{n-1}(x)-1\}.$$

Because $t'_n(x)=0$ for any $x$ outside
$\bigcup_{U\in\mathcal{U}_n(P_0)}U^1\setminus U$, we observe that
\begin{equation*}
t_{n}(x) =
\begin{cases}
t'_n(x) & \text{ if }x\in\bigcup_{U\in\mathcal{U}_n(P_0)}U^1\setminus U\\
0 & \text{ if }x\in \AA_{n-1}(P_0)\setminus\big(\bigcup_{U\in\mathcal{U}_n(P_0)}U^1\setminus U\big)\\
t_{n-1}(x)-1 & \text{ if }x\in\BB_{n-1}(P_0)\setminus\big(
\bigcup_{U\in\mathcal{U}_n(P_0)}U^1\setminus U\big)
\end{cases},
\end{equation*}
whenever $\,x\in\Delta_n(P_0)$. This finishes the second case.

In both cases, $\Delta_n(P_0)$ and $t_n$ are defined  and we set
$$\AA_{n}(P_0)= \{x\in\Delta_{n}(P_0)\,\|\, t_{n}(x) = 0\}$$
$$\mbox{and}$$
$$\BB_{n}(P_0)=
\{x\in\Delta_{n}(P_0)\,\|\, t_{n}(x) > 0\}.$$

\begin{Remark}\label{rem1} As $t_n=\max\{t'_n,t_{n-1}-1\}$ $\forall\,n\ge
1$, we get $t_n\ge t_{n-1}-1\ge t_{n-2}-2\ge...\ge t_{n-j}-j$. So,
$$t_{n}\ge t_k-(n-k)\mbox{ }\forall n\ge k\ge 0.$$
\end{Remark}

\begin{Remark}\label{rem2}
Note that $\dist(U^1,M\setminus P_0)>\delta_0 \sigma^{n/2}$,
$\forall\, U\in\mathcal{U}_n(P_0)$ and $n\ge1$. In particular,
$$U^1\subset P_0,\,\forall\, U\in\mathcal{U}_n(P_0).$$
\end{Remark}

To check the remark~\ref{rem2}, first observe that
$\dist(\AA_{n-1}(P_0),M\setminus P_0)\ge\delta_0 \sigma^{(n-1)/2}$
(because $t_{n-1}(x)=0$ implies by remark~\ref{rem1} that $t_0(x)\le
n-1$). On the other hand, $U^1$ is a pre-image (at a hyperbolic
time) of some set $B_{\delta_0}(P)$ ($P\in\cp$) with diameter
smaller than $3\delta_0$, and so the diameter of $U^1$ is smaller
than $3\delta_0\sigma^{n/2}$. Since
$U^1\cap\AA_{n-1}(P_0)\ne\emptyset$, we get $\dist(U^1,M\setminus
P_0)\ge\delta_0 \sigma^{(n-1)/2}-3\delta_0 \sigma^{n/2}=\delta_0
\sigma^{n/2}(1/\sigma^{1/2}-3)>\delta_0 \sigma^{n/2}$.

At this point we have completely described the inductive
construction restricted to $P_0$ of the sets  $\AA_n(P_0)$,
$\BB_n(P_0)$ and $\Delta_n(P_0)$. Proceeding in the same way to the
other elements of $\cp$ we obtain
$$\AA_n=\bigcup_{P\in\cp}\AA_n(P),$$
$$\BB_n=\bigcup_{P\in\cp}\BB_n(P)$$
$$\mbox{and}$$
$$\Delta_n=\bigcup_{P\in\cp}\Delta_n(P).$$

Let $M^*=\{x\in M\,\|\, f^n(x)\in\Delta_0\,\forall\,n\ge0\}$. As $f$
is a diffeomorphism  Lebesgue almost everywhere, we have
$\leb(M\setminus M^*)=0$. Finally let us define a {\em time
function}
$$R:M^*\setminus\bigcap_n\Delta_n\to\NN$$ setting $R(x)=n$ when
$x\in\Delta_{n-1}\setminus\Delta_n$.

We may think of the set $\AA_{n-1}$
($\{t_{n-1}=0\}\cap\Delta_{n-1}$) as the set of {\em allowed points}
of the $n$-th step of the construction, that is, the set of points
$x\in M$ that can be used in the step n, if $n$ is a hyperbolic time
to $x$, as the ``center'' of a new element of the partition (a
component of $\{R=n\}$). On the other hand, the set $\BB_{n-1}$
($\{t_{n-1}>0\}\cap\Delta_{n-1}$) is the {\em forbidden points} of
the $n$-th step, i.e., if $x\in \BB_{n-1}$ then, though $x$ does not
belong to a constructed element ($x\notin\{R\le {n-1}\}$), we can
not use $x$ to be the ``center'' of a new element, even if $x\in
H_{n}$. In fact, $\BB_{n-1}$ is the intersection of $\Delta_{n-1}$
with an union of ``protection collars''  associated to the
components of $\{R\le {n-1}\}$ and we can not use the points of
$\BB_{n-1}$ in the step $n$ because at this time the pre-image of
the elements of the triangulation associated to some hyperbolic
pre-ball $V_{n}(x)$ (with $x\in H_{n}$) may be big compared with
these {\em collars} and, probably, there will be overlaps of this
pre-image with some element constructed already.

\begin{Remark}\label{preball}It is important to emphasize that we
only constructed a new element of the partition as a pre-image,
associated to a hyperbolic time, of some element of the
triangulation. That is, if $U$ is a component of $\{R=n\}$, for some
$n$, then $U\subset V_n(x)$ for some $x\in\AA_{n-1}\cap H_n$. As a
consequence, $U$ is diffeomorphically mapped, with bounded
distortion and good properties of expansion (see
proposition~\ref{Prop3}), to some element of the triangulation. This
control of distortion and expansion is fundamental to obtain the
properties required in the piecewise uniformly expanding Markovian
map associated to the Markov structure.
\end{Remark}

\subsection{Collars and rings}To make precise the meaning of
``collar'', we introduce some additional notation. Given
$P_0\in\mathcal P$ we define, for any $n\ge 0$, the {\em $n$-th
former collar of $P_0$}, or the {\em $(0,n)$-collar of $P_0$} as the
set
\begin{equation}\label{former1}
C_{0,n}(P_0)=\{x\in P_0\,\|\,\dist(x,\partial P_0)<\delta_0
\sigma^{n/2}\},
\end{equation}
that is,
\begin{equation}\label{former2}
C_{0,n}(P_0)=\{t_0\ge n+1\}\cap P_0.
\end{equation}

For any element $U\in\mathcal U_k(P_0)$ (a connected component of
$\{R=k\}$) and each $n\ge k$, we define the {\em $(k,n)$-th collar
of $P_0$ associated to $U$} as
$$C_{k,n,U}(P_0)=\bigcup_{j=n-k+1}^\infty U^1(j).$$

A {\em collar} of the $n$-th step of the construction is any
$(k,n)$-th collar, with $k\le n$. The collection of all collars of
the $n$-th step associated to $P_0$ is denoted by $\beta_n(P_0)$,
i.e.,
$$\beta_n(P_0)=\{C_{0,n}(P_0)\}\cup\{C_{k,n,U}(P_0)\,\|\,
U\in\mathcal U_k(P_0)\mbox{ and }1\le k\le
n\}.$$

Finally, we define $C'_{k,n,U}(P_0)$, the {\em external ring} of the
collar $C_{k,n,U}$ $(P_0)$, as
\begin{equation}\label{ring1}
C'_{k,n,U}(P_0)=\underbrace{C_{k,n,U}(P_0)\setminus
C_{k,n+1,U}(P_0)}_{U^1(n-k+1)}
\end{equation}
and the {\em external ring} of the former collar $C_{0,n}(P_0)$ as
\begin{equation}\label{ring2}
C'_{0,n}(P_0)=\underbrace{C_{0,n}(P_0)\setminus
C_{0,n+1}(P_0)}_{C_{0,n}(P_0)\cap \{t_0=n+1\}}.
\end{equation}

We claim that
\begin{equation}\label{collar}
\bigcup_{Q\in\beta_{n}}Q=\{t_{n}\ge
1\}\hspace{0.5cm}\mbox{and}\hspace{0.4cm}\bigcup_{Q\in\beta_{n}}Q\setminus
Q'=\{t_{n}> 1\}\hspace{0.3cm}\forall n\in\NN,
\end{equation}
where $\beta_{n}=\bigcup_{P\in\cp}\beta_n(P)$ and $Q'$ is denoting
the external ring of the collar $Q$. Indeed, if $x\in
C_{0,n}(P_0)=\{t_0\ge n+1\}\cap P_0$, it follows form
Remark~\ref{rem1} that $t_n(x)\ge 1$. On the other hand, if $x\in
C_{k,n,U}(P_0)$, then $x\in U^1(j)$ for some $j\ge n-k+1$. Hence,
$t_n(x)\ge t_k(x)-(n-k)\ge t'_k(x)-(n-k)\ge j-(n-k)\ge 1$. In the
same way, if $x\in C_{0,n}(P_0)\setminus C'_{0,n}(P_0)$ or
$C_{k,n,U}(P_0)\setminus C'_{k,n,U}(P_0)$ we get $t_n(x)>1$. So,
$\bigcup_{Q\in\beta_{n}}Q\subset\{t_{n}\ge 1\}$ and
$\bigcup_{Q\in\beta_{n}}Q\setminus Q'\subset\{t_{n}> 1\}$. Now, let
us assume that $x\in P_0$ and $t_n(x)=j\ge 1$. Note that
$t_n=\max\{t'_n, t_{n-1}-1\}$ $=$ $\max\{t'_n, \max\{t'_{n-1},
t_{n-2}-1\}-1\}$ $=\max\{t'_n, t'_{n-1}-1, t_{n-2}-2\}=$
$\max\{t'_n, t'_{n-1}-1, t'_{n-2}-2, t_{n-3}-3\}=\ldots$. Thus,
\begin{equation}\label{equa1}
t_n=\max\{t'_n, t'_{n-1}-1, t'_{n-2}-2, \ldots, t'_{1}-(n-1),
t_0-n\}
\end{equation}
and  $t_n(x)=j\ge 1$ implies that $t_0(x)=j+n$ or $t'_k(x)=j+(n-k)$,
for some $1\le k\le n$. In the first case we get $x\in C_{0,n}(P_0)$
(and $x\in C'_{0,n}(P_0)$ when $j=1$). In the second one, by the
definition of $t'_k$, there exists $U\in\mathcal U_k(P_0)$ such that
$x\in U^1(n-k+j)$. So, $x\in C_{k,n,U}(P_0)$ (and $x\in
C'_{k,n,U}(P_0)$ only if $j=1$). In any case, we obtain $\{t_{n}\ge
1\}\subset\bigcup_{Q\in\beta_{n}}Q$ (and also $\{t_{n}>
1\}\subset\bigcup_{Q\in\beta_{n}}Q\setminus Q'$), thus proving the
claim.

The next lemma is important for preventing the overlaps on the sets
of the partition. Indeed, in the $n$-th step of the induction,
associated to each element $V\in\cu_k(P_0)$ already constructed ($V$
is a connected component of $\{R=k\}$, for some $0<k<n$), there is a
collar $Q=C_{k,n,V}(P_0)\in\beta_n(P_0)$ around it. By the lemma
below, the new components cannot intersect ``too much'' the collar
$Q$. More precisely, if $U\in\mathcal{U}_n(P_0)$ (i.e.,$U$ is a
connected component of $\{R=n\}$) then the possible intersection of
$U$ with $Q$ will be necessarily within the external ring
$Q'=C'_{k,n,V}(P_0)$. As a consequence, according to
Corollary~\ref{disjointed}, there will be no overlaps.

\begin{Lemma}\label{l.claim} If
$U\in\mathcal{U}_n(P_0)$ then $U\cap\{t_{n-1}>1\}=\emptyset$.
\end{Lemma}

\dem  Let us assume that $U\cap\{t_{n-1}>1\}\ne\emptyset$ for some
$U\in\mathcal{U}_n(P_0)$. The set $U$ can be written as
$U=(f^n|_{V_n(x_1)})^{-1}(P_1)$, where $x_1$ is some point in
$\AA_{n-1}\cap H_n$ and $P_1$ is an element of the triangulation
$\cp$. By (\ref{collar}), we have
\begin{equation}\label{222}
U\cap (Q\setminus Q')\ne\emptyset\,\mbox{ for some collar
}Q\in\beta_{n-1}(P_0),
\end{equation}
where  $Q'$ is its external ring. Firstly let us verify that $Q$ is
not the $(n-1)$-th former collar, i.e, $Q\ne C_{0,n-1}(P_0)$. For
this, let us suppose by contradiction that $Q=C_{0,n-1}(P_0)$. It
follows from Remark~\ref{rem2} that $\dist(U,\partial P_0)>\delta_0
\sigma^{n/2}$ and, as a consequence, $t_0(x)\le n$ $\forall\,x\in
U$. So, $U\cap(Q\setminus Q')=U\cap \big(C_{0,n-1}(P_0)\setminus
C'_{0,n-1}(P_0)\big)=U\cap \{t_0>n\}\cap P_0=\emptyset$, thus
contradicting (\ref{222}).

As $Q$ is not the $(n-1)$-th  former collar,  there exist $1\le k\le
n-1$ and $V\in\mathcal{U}_k(P_0)$ such that
$Q=C_{k,n-1,V}(P_0)\in\beta_{n-1}(P_0)$. On its turn, $V$ can be
written as $(f^k|_{V_k(x_2)})^{-1}(P_2)$, where $x_2$ is some point
of $\AA_{k-1}\cap H_k$ and $P_2\in\cp$.

By construction, $f^k$ sends diffeomorphically $V$ onto $P_2$,
$V\cup Q$ onto the topological ball $B_{\delta_0
\sigma^{(n-1-k)/2}}(P_2)$ and $Q'$ onto the ``ring''
$I^1_{n-k}(P_2)=B_{\delta_0 \sigma^{(n-1-k)/2}}(P_2)\setminus
B_{\delta_0 \sigma^{(n-k)/2}}(P_2)$. Thus, $Q'$ splits $M$ into two
connected components, that is, taking $W=V\cup(Q\setminus Q')$ and
$L=M\setminus(V\cup Q)$, we get $M\setminus Q'=W\cup L$ with $W\cap
L=\emptyset$. Note that $W\supset V$ and $L\supset\AA_{n-1}$
(because $V\cap\AA_{j}\subset V\cap\Delta_j=\emptyset$
$\forall\,j\ge k$ and, as $t_k(x)\ge t'_k(x)\ge n-k$ $\forall\,x\in
Q$, it follows from remark~\ref{rem1} that $t_{n-1}(x)\ge
n-k-(n-1-k)=1>0$ $\forall\,x\in Q$).

We know that $U\cap\AA_{n-1}\ne\emptyset$, because, by construction,
$x_1\in\AA_{n-1}\cap H_n$. Hence, $U$ intersects $L$. On the other
hand, we are assuming that $U\cap W\supset U\cap(Q\setminus
Q')\ne\emptyset$. So, it follows from the connectivity of $U$ that
it not only intersects $Q'$ but also intersects both connected
components of the boundary of $Q'$ (recall that $f^n$ sends
diffeomorphically $\partial Q'$ onto $\partial
I^1_{n-k}(P_2)=\partial B_{\delta_0
\sigma^{(n-1-k)/2}}(P_2)\cup\partial B_{\delta_0
\sigma^{(n-k)/2}}(P_2)$).

Therefore, we can pick two points in $U$, $q_1$ and $q_2$, in
distinct components of the boundary of $Q'$. We assume that
$f^k(q_1)\in\partial B_{\delta_0 \sigma^{(n-1-k)/2}}(P_2)$ and
$f^k(q_2)\in\partial B_{\delta_0 \sigma^{(n-k)/2}}(P_2)$. As $U$ is
the pre-image of $P_1$ in a hyperbolic time, we get from
Proposition~\ref{Prop3} that
 \begin{equation}\label{e.zq1}
 \dist(f^{k}(q_1),f^{k}(q_2))\le
\sigma^{(n-k)/2}\dist(f^{n}(q_1),f^{n}(q_2))<\delta_0
\sigma^{(n-k)/2}.
 \end{equation}
On the other hand,
 \begin{eqnarray*}
\dist(f^{k}(q_1),f^{k}(q_2))&\ge&
\delta_{0}\,\sigma^{(n-1-k)/2}-\delta_{0}\,\sigma^{(n-k)/2}\\&=&
\delta_{0} \sigma^{(n-k)/2}(\sigma^{-1/2}-1),
 \end{eqnarray*}
which combined with (\ref{e.zq1})  gives $\sigma>\frac{1}{4}$. This
contradicts our choice of $\sigma$ (see (\ref{sigma})). \cqd

\begin{Corollary}\label{disjointed} If $U$ and $V$ are distinct
components  of $\{R=n\}$ and $\{R=m\}$ then $U\cap V=\emptyset$
\end{Corollary}
\dem Suppose that $U=(f^n|_{V_n(x_1)})^{-1}(P_1)$ and
$V=(f^m|_{V_m(x_2)})^{-1}(P_2)$, with $P_1$  and $P_2\in\cp$. By
Remark~\ref{rem2} we can assume that both $U$ and $V$ are contained
in some $P_0\in\cp$. As we claimed during the presentation of the
algorithm, if $n=m$ we have $U\cap V=\emptyset$. So, let us assume
that $m<n$. Let $Q=C_{m,n-1,V}(P_0)$ be the $(n-1)$-th collar
associated to $V$ and  $Q'=C'_{m,n-1,V}(P_0)$ its external ring.

By construction, $f^m$ sends diffeomorphically $V$ onto the
topological ball $P_2$ and the set $T=Q\setminus Q'$
diffeomorphically onto the ``ring'' $B_{\delta_0
\sigma^{(n-m)/2}}(P_2)\setminus
P_2=\bigcup_{k=n-m+1}^{\infty}I^1_{k}(P_2)$. Thus, $T$ splits $M$
into two connected components. Indeed, taking $L=M\setminus(V\cup
T)$, we get $M\setminus T=V\cup L$ with $V\cap L=\emptyset$.

It is easy to see that $T\subset\{t_{n-1}>1\}$. In fact, as
$t_m(x)\ge t'_{m}(x)\ge n-m+1$ $\forall\,x\in T$ ($=Q\setminus
Q'=V^1((n-1)-m+2)\cup V^1((n-1)-m+3)\cup\ldots$), it follows from
Remark~\ref{rem1} that $t_{n-1}(x)\ge n-m+1-(n-1-m)=2$
$\forall\,x\in T$.

Note that $L\supset\AA_{n-1}$ (because $V\cap\AA_{j}\subset
V\cap\Delta_j=\emptyset$ $\forall\,j\ge m$ and $\AA_{n-1}\cap
T\subset \AA_{n-1}\cap\{t_{n-1}\ne 0\}=\emptyset$). We know that
$U\cap\AA_{n-1}\ne\emptyset$, because by construction
$x_1\in\AA_{n-1}\cap H_n$. Hence, $U$ intersects $L$. If $U\cap
V\ne\emptyset$, it follows from the connectivity of $U$ that $U\cap
T\ne\emptyset$. But this contradicts Lemma~\ref{l.claim}, as
$T\subset\{t_{n-1}>1\}$. \cqd

The proposition below says that the sum of the Lebesgue measure of
all forbidden sets is finite. This finiteness is fundamental to
assure that almost every point in $M$  belongs to a constructed
element, i.e., $\leb(M\setminus\bigcup_n\{R=n\})=0$.

\begin{Proposition}\label{PropB}
$\sum_{n=0}^\infty \leb\,\BB_{n}<\infty$.
\end{Proposition}
 \dem

First of all, we compare the Lebesgue measure of the $(k,n)$-th
collar $C_{k,n,U}(P_0)$ with the measure of $U\in\cu_k(P_0)$, where
$P_0$ is an element of $\cp$ and $1\le k\le n$. For this, recall
that there must be some hyperbolic pre-ball $V_k(x_0)$ (containing
$U$ and $C_{k,n,U}(P_0)$) and an element $P_1$ of $\cp$ such that
$f^k|_{V_k(x_0)}$ maps $U$ diffeomorphically onto $P_1$ and
$C_{k,n,U}(P_0)$ onto $\bigcup_{j=n-k+1}^\infty I_j^1(P_1)$, both
with the distortion  bounded  by an universal constant
$C_1=\exp(\rho\,\mbox{diameter}(M))$ (see proposition~\ref{Prop3}).
Thus, it is sufficient to compare the measures of
$\bigcup_{j={n-k+1}}^\infty I_j^1(P_1)$ and $P_1$. As $P_1\in\cp$ is
the interior of an element of a triangulation of $M$, the boundary
of $P_1$ is a finite union of smooth codimension one submanifolds.
So, for $r>0$ small we have $\leb(B_{r}(P_1)\setminus P_1)\approx
r$. As $\cp$ is a finite collection, one can find $\xi>0$ such that
\begin{equation}\label{comp1}
\frac{\leb(B_{r}(P)\setminus P)}{\leb\,P}<\xi
r\hspace{0.3cm}\mbox{for all } P\in\cp\mbox{ and }\,r>0.
\end{equation}
Hence,
$$\frac{\leb\,C_{k,n,U}(P_0)}{\leb\,U}\le
C_1\frac{\leb(\bigcup_{j={n-k+1}}^\infty I_j^1(P_1))}{\leb\,P_1}=$$
$$=C_1\frac{\leb(B_{\delta_0\sigma^{n-k}}(P_1)\setminus P_1)}{\leb\,P_1}
< C_1 \xi \delta_0\sigma^{n-k}.$$ Thus, setting $C_2=C_1 \xi
\delta_0$, we have
\begin{equation}\label{comp2}
\leb\,C_{k,n,U}(P_0)<C_2 \sigma^{n-k} \leb\,U
\end{equation}
for all $1\le k\le n$, $U\in\cu_n(P_0)$ and $P_0\in\cp$.

As in (\ref{comp1}), one can find some $\varsigma>0$ such that
\begin{equation}\label{comp3}
\frac{\leb\,C_{0,n}(P)}{\leb\,
P}<\varsigma\sigma^{n}\hspace{0.3cm}\mbox{for all }P\in\cp\mbox{ and
}\,n\ge 0.
\end{equation}

Using (\ref{comp2}) and (\ref{comp3}), it follows that
$$\sum_{Q\in\beta_n(P_0)}\leb\,Q=\leb\,C_{0,n}(P_0)+
\sum_{k=1}^n\sum_{U\in\cu_k(P_0)}\leb\,C_{k,n,U}(P_0)$$
$$<\varsigma\sigma^{n}\leb\,P_0+\sum_{k=1}^n\sum_{U\in\cu_k(P_0)}C_2 \sigma^{n-k}
\leb\,U=$$
$$=\varsigma\sigma^{n}\leb\,P_0+C_2\sum_{k=1}^n\sigma^{n-k}\leb(\{R=k\}\cap P_0).$$
Thus,
$$\sum_{Q\in\beta_n}\leb\,Q=\sum_{P_0\in\cp}\sum_{Q\in\beta_n(P_0)}\leb\,Q
<$$
$$<\varsigma\sigma^{n}\underbrace{\leb\,M}_1+C_2\sum_{k=1}^n\sigma^{n-k}\leb\{R=k\},
\hspace{0.4cm}\forall n\ge1.$$

As
$$\sum_{n=1}^{\infty}\sum_{k=1}^n\sigma^{n-k}\leb\{R=k\}\le
\leb\{R=1\}+$$ $$+\big(\sigma \leb\{R=1\}+\leb\{R=2\}\big)+$$
$$+\big(\sigma^2 \leb\{R=1\}+\sigma \leb\{R=2\}+\leb\{R=3\}\big)+\ldots=$$
$$=\sum_{j=0}^{\infty}\sigma^j\,\leb\{R=1\}+\sum_{j=0}^{\infty}\sigma^j\,
\leb\{R=2\}+\ldots=$$ $$=\frac{1}{1-\sigma}\sum_{n=1}^{\infty}
\leb\{R=n\}\le\frac{1}{1-\sigma}\underbrace{\leb\,M}_{1}\,,$$ we get
$$\sum_{n=0}^{\infty}\sum_{Q\in\beta_n}\leb\,Q=\varsigma+
\sum_{n=1}^{\infty}\sum_{Q\in\beta_n(P_0)}\leb\,Q\le\varsigma+
\varsigma\frac{\sigma}{1-\sigma}+\frac{C_2}{1-\sigma}<\infty.$$ This
completes the proof, because
$\BB_n=\bigcup_{Q\in\beta_n}Q\cap\Delta_n$, and so,
$$\sum_{n=0}^\infty
\leb\,\BB_{n}<\frac{\varsigma+C_2}{1-\sigma}<\infty.$$

\cqd

Now, we are going to prove that this algorithm does indeed produce a
partition Lebesgue mod 0 of $M$. As Corollary~\ref{disjointed}
assures that the elements constructed by this algorithm are mutually
disjoint, we have only to check that $\leb(\bigcap_n\Delta_n)=0$
($\Delta_0\supset\Delta_1\supset...$ is a nested sequence).

To prove this, first observe that, as $\sum_{n=0}^\infty
\leb\,\BB_{n}<\infty$, it follows from Borel-Cantelli lemma that
Lebesgue almost every point in $M$ belongs only to  finitely many
${\BB_j}'s$, that is, defining the function $b:M\to\NN$ by
\begin{equation}\label{b}
    b(x) =
    \begin{cases}
    \#\{n\in\NN\,\|\,x\in\BB_{n-1}\} & \text{ if } x\in\Delta_0\\
    +\infty & \text{ if }x\notin\Delta_0\\
    \end{cases},
\end{equation}
we have $b(x)<+\infty$ for Lebesgue almost all $x\in M$ (because
$\int b\,d Leb=\sum_n\leb\,\BB_n$). Moreover, as any Lebesgue
generic point has infinity many hyperbolic times, one can find for
almost every $x\in M^*$ the first time
$n>\max\{n\in\NN\,\|\,x\in\BB_{n-1}\}$ such that $x\in H_n$. In this
case, if $x$ still belongs to $\Delta_{n-1}$ we have
$x\in\Delta_{n-1}\setminus\Delta_{n}$. In other words, almost every
point in $M^*$ (and consequently, in $M$) belongs to some
$\Delta_{n-1}\setminus\Delta_{n}$. So, $\leb(\bigcap_n\Delta_n)=0$.

\subsection*{Markov Structures}\label{Markov Structures}

Let $F:\Delta^*\to M$ be given by $F(x)=f^{R(x)}(x)$, where
$\Delta^*=M^*\setminus\bigcap_n\Delta_n$. By construction, the map
$F$ is a {\em piecewise uniformly expanding Markovian map\/} (see
remark~\ref{preball}) and we conclude the proof of the first part of
the Main Theorem.

\section{Integrability of The Time Function}

In this section we finish the proof of the Main Theorem showing the
integrability of the time function of the Markov Structure
constructed in the previous section. We will use the objects and
notation of the {\em Partitioning Algorithm}. Let $\nu$ be one of
the $F$-invariant measure given by the Markov Structure and assume
by contradiction that $\sum_{P\in\mathcal{P}}R|_P\,\nu(P)=\infty$.
It follows from Birkhoff's Ergodic Theorem that
$$
\frac{1}{n}\sum_{j=0}^{n-1}R\circ F^j(x)\longrightarrow\int
R\,d\nu=\sum_{P\in\mathcal{P}}R|_P\,\nu(P)=\infty,
$$
for $\nu$-almost every point $x\in M$.

As the density of $\nu$ is uniformly bounded from above, it follows
from proposition~\ref{PropB} that
$\sum_{n=0}^{\infty}\nu(\BB_n)<\infty$. Hence, for $\nu$-generic
points $x\in M$ we get
$$\frac{1}{n}\sum_{j=0}^{n-1} b\circ F^j(p)\longrightarrow\int
b\,d\nu=\sum_n\nu(\BB_n)<\infty$$

Let $M^{**}=\{x\in M^*\,\|\,F^j(x)\in Dom(F)\,\forall\,j\ge 0\}$. It
is clear that $\nu(M\setminus M^{**})=0$. Let $x\in M^{**}$. Set,
for every $i\in\NN$,  $j_i=j_i(x)=\sum_{j=0}^{i-1}R\circ F^j(x)$,
that is, $F^i(x)=f^{j_i}(x)$. Given  $j\in\NN$, there exists a
unique integer $s=s(j)\ge0$ such that $j_s\le j < j_{s+1}$. Let us
set $\mathcal{I}_0=\{j_1,j_2,,j_3,...\}$ and suppose that $x\in
H_j$. In this case, $F^{s}(x)\in H_{m}$, where $m=j-j_s$. By
construction, if $F^{s}(x)\notin\{R=m\}$ (i.e. $m<R(F^s(x))$), then
$F^{s}(x)\in\BB_{m-1}$. On the other hand, if $F^{s}(x)\in\{R=m\}$,
then $m=j_{s+1}-j_s$, which implies  $j=j_{s+1}\in\mathcal{I}_0$.

As the number of integers $\ell$ between $j_i$ and $j_{i+1}$ such that
$x\in H_\ell$ is bounded by the number of integers $m$  such that
$F^i(x)\in \BB_{m-1}$, we have
$$\#\{\ell\in\{j_i+1,...,j_{i+1}-1\}\,\|\,x\in H_\ell\}\le b(F^{i}(x)).$$

Thus, for each $n\in\NN$ we can write
$$\#\{j\le n\,\|\,x\in H_j\}\le\#\{j\le n\,\|\,j\in\mathcal{I}_0\}+
\sum_{i=0}^{s(n)} b(F^i(x))\le$$

$$\le s(n)+\sum_{i=0}^{s(n)} b(F^i(x))$$

Therefore,
\begin{equation}\label{eqS1}
\frac{1}{n}\#\{1\le j\le n\,\|\,x\in H_j\}\le
\frac{s(n)}{n}\bigg(1+\frac{1}{s(n)} \sum_{j=0}^{s(n)} b\circ F^j(x)\bigg)
\end{equation}
By construction,  if $s(n)=i$, that is, $j_i\le n <j_{i+1}$, then
\begin{equation}\label{eqS2}
\frac{j_i}{i}\le\frac{n}{s(n)}<\frac{j_{i+1}}{i+1}\bigg(1+\frac{1}{i}\bigg)
\end{equation}
As
$\frac{j_i}{s(j_i)}=\frac{j_i}{i}=\frac{1}{i}\sum_{j=0}^{i-1}R\circ
F^j(x)\to\infty$, it follows from equations (\ref{eqS1})~and
(\ref{eqS2}) that
$$\lim_{n\to\infty}\frac{1}{n}\#\{1\le j\le n\,\|\,x\in
H_j\}=0.$$ But this is an absurd as one can see in
Remark~\ref{remarkTHETA}. So, we necessarily have the time function
$R$ $\nu$-integrable and so, we conclude the proof of the {\em Main
Theorem}.

To prove theorems~\ref{SRB}, \ref{weaklyXnue} and \ref{HTSRB} let
$\{\nu_1,...,\nu_s\}$ be the set of ergodic absolutely continuous
invariant measures given by Theorem~\ref{MST}.

\subsection*{Proof of Theorem~\ref{SRB}}
 As we observed before it is
straightforward to check that as the time function $R$ is $\nu_i$-integrable,
then
$$
\mu_i =\sum_{j=0}^{\infty}f_{\ast}^j\left(\nu_i\vert \{R>j\}\right)
$$
is an ergodic absolutely continuous $f$-invariant measure. Here
$\nu_i\vert \{R>j\}$  denotes the measure given by $\nu_i\vert
\{R>j\}(A)=\nu_i(A\cap\{R>j\})$, and $f_*^j$ denotes the
push-forward of the measure by $f^j$. Moreover, it follows from
Theorem~\ref{MST} that almost every point in $M$ belongs to the
basin of one of the $\nu_i$. Therefore, almost every point $x\in M$
also belongs to the basin of one of the measures $\mu_1,...,\mu_s$.
\begin{flushright}\mbox{$\square$}\end{flushright}

\begin{Remark}
In the beginning of the Partitioning Algorithm we could have
replaced, if necessary, $f$ by some iterated  $f^{n_0}$ and so, the
measures $\mu_1,...,\mu_s$ might  be only $f^{n_0}$ invariant. In
such case we consider the induced $f$-invariant measures
$\mu_1',...,\mu_s'$ given by
$$\mu_i'(A)=\mu_i(A)+\mu_i(f^{-1}A)+...+\mu_i(f^{-(n_0-1)}A)$$
\end{Remark}

\subsection{Proof of Theorem~\ref{weaklyXnue}}
As any non-uniformly expanding map satisfies the hypothesis of
Theorem~\ref{SRB} we have only to check that every non-flat map with
the slow recurrence satisfying equation~(\ref{e3}) for almost
every$x\in M$ also satisfies equation~(\ref{e2}). This result is a
consequence of the integrability of $x\mapsto\log\|(Df(x))^{-1}\|$
with respect to any absolutely continuous invariant measure ( see
Lemma~\ref{PDL} below ). Indeed, as
$\log\|(Df)^{-1}\|^{-1}=-\log\|(Df)^{-1}\|$, if
$\,\log\|(Df)^{-1}\|\in L^{1}(\mu_i)$ it follows from Birkhoff's
Ergodic Theorem that we have for $\mu_i$ almost all $x\in M$
$$\liminf_{n\to\infty}\frac{1}{n}\sum_{j=0}^{n-1}\log\|(Df(f^n(x)))^{-1}\|^{-1}=$$
$$=\limsup_{n\to\infty}\frac{1}{n}\sum_{j=0}^{n-1}\log\|(Df(f^n(x)))^{-1}\|^{-1}.$$
As Lebesgue almost every point $x\in M$ belongs to the basin of some
$\mu_i$, we conclude that the equality above of $\liminf$ and
$\limsup$ holds for Lebesgue almost all $x\in M$.
\begin{flushright}\mbox{$\square$}\end{flushright}

\begin{Lemma}\label{PDL}
If $f:M\to M$ is a $C^2$ endomorphism, then $\log\|(Df)^{-1}\|$ is
integrable with respect to any absolutely continuous invariant
measure.
\end{Lemma}
\dem
In \cite{Liu} (see remark 1.2 of \cite{Liu}) it was proved that
$$-\infty<\int\log|\det Df| d\mu<+\infty$$ with respect to any
absolutely continuous invariant measure $\mu$. Let us show how this
result implies the integrability of $\log\|(Df)^{-1}\|$. It is easy
to see that $|\lambda|\le\|A\|$ for any eigenvalue $\lambda$ of a
matrix $A\in\RR^{m\times m}$. Thus, $|\det A|\le\|A\|^m$ and so
$$-\infty<\frac{1}{dim(M)}\int\log|\det
Df|d\mu\le\int\log\|Df\|d\mu<+\infty.$$ We know that $\adj(A)\,
A=\det(A)\, I$ $\forall A\in\RR^{n\times n}$, where $\adj\, A$ is
the (classical) adjoint matrix, the transposed matrix of cofactors
\cite{HJ}, and $I$ is the identity matrix. So, $\log|\det
Df|-\log\|Df\|\le\log\|\adj{Df}\|$ and so
$-\infty<\int\log\|\adj{Df}\|d\mu$. Moreover, as $\|\adj\,Df\|$ is
bounded, we get
$$-\infty<\int\log\|\adj\,Df\|d\mu<+\infty.$$ Finally, from
$(Df)^{-1}=\frac{1}{\det Df}\,adj\,Df$ we get
$$-\infty<\int\log\|(Df)^{-1}\|d\mu=$$
$$=\int\log\|\adj\,Df\|d\mu-\int\log|\det Df|d\mu<+\infty.$$ \cqd

\subsection{Proof of Theorem~\ref{HTSRB}}
Suppose that there exists $0<\sigma<1$ and $\delta>0$ such that
Lebesgue almost every point $x\in M$ has a
$(\sigma,\delta)$-hyperbolic time. So, we have $\leb(U)=0$, where
$U=M\setminus\bigcup_j H_j(\sigma,\delta)$. Given a Lebesgue generic
point $x\in M$ we have $f^j(x)\notin U$ $\forall j$. Let $n_1$ be
such that $x\in H_{n_1}(\sigma,\delta)$ and $m_1$ such that
$f^{n_1}(x)\in H_{m_1}(\sigma,\delta)$. By the properties of
hyperbolic times we get $x\in H_{n_2}(\sigma,\delta)$, where
$n_2=n_1+m_1>n_1$. Therefore, by induction, one concludes that
Lebesgue almost every point in $M$ has infinitely many
$(\sigma,\delta)$-hyperbolic times. Remembering that $x\in
H_j(\sigma,\delta)$ implies $\prod_{k=0}^{j-1}\|(Df\circ
f^k(x))^{-1}\|\le{\sigma}^j$, we have, for generic $x\in M$,
$\frac{1}{j}\sum_{k=0}^{j-1}\log\|(Df(f^k(x)))^{-1}\|\le\log\sigma<0$
for infinitely many integers $j>0$. Therefore,
$$\liminf_{j\to\infty}\frac{1}{j}\sum_{k=0}^{j-1}\log\|(Df(f^k(x)))^{-1}\|\le\log\sigma<0$$
for Lebesgue almost all $x\in M$. As the critical set is empty, the
map $f$ satisfies the hypothesis of the theorem~\ref{SRB}. So, by
the theorem~\ref{weaklyXnue}, $f$ is a non-uniformly expanding map.
\begin{flushright}\mbox{$\square$}\end{flushright}

\begin{Remark}
Using the same argument as in the proof of Theorem~\ref{HTSRB}, it
is easy to conclude that we can change, in the hypothesis of
Theorem~\ref{HTSRB}, the existence of a hyperbolic time by the
existence of a time with good average of expansion. That is, if
$f:M\to M$ is $C^2$ local diffeomorphism and $\exists\,\lambda>0$
such that for Lebesgue almost every point $x\in M$ one can find at
least one $n_0=n_0(x)\in\NN$ such that
$$
\frac{1}{n_0}\sum_{j=0}^{n_0-1}\log(\|(Df(f^j(x)))^{-1}\|^{-1})\ge\lambda>0,
$$
then $f$ is a non-uniformly expanding map.
\end{Remark}

\end{document}